\documentclass[11pt,reqno]{amsart}
\usepackage{amsmath,amssymb,amsthm}
\usepackage{mathrsfs}
\usepackage{enumitem}
\usepackage[margin=1in]{geometry}
\usepackage[hidelinks]{hyperref}
\usepackage{comment}
\usepackage{xcolor}
\DeclareMathOperator{\supp}{supp}

\newcommand{\R}{\mathbb{R}}
\newcommand{\F}{\mathcal{F}}
\newcommand{\Fi}{\mathcal{F}^{-1}}

\theoremstyle{plain}
\newtheorem{theorem}{Theorem}[section]

\theoremstyle{definition}

\newtheorem{conjecture}{Conjecture}
\theoremstyle{remark}

\newcommand{\step}[1]{\medskip\par\noindent\textit{#1}\quad}
\title[The fundamental solution of the logarithmic Laplacian is a function]
{The fundamental solution of the logarithmic Laplacian is a function in dimensions 1 and 2}

\author{David Lee}
\address{School of Mathematical and Physical Sciences,
University of Technology Sydney (UTS)}
\email{davidchanwoo.lee@uts.edu.au}

\author{Samuel Killick}
\address{School of Mathematical and Physical Sciences,
University of Technology Sydney (UTS)}
\email{samuel.killick@uts.edu.au}

\date{\today}


\begin{document}

\begin{abstract}
We show that the fundamental solution of the logarithmic Laplacian in dimensions $1$ and $2$ is a locally integrable function. Together with earlier work of the first author, this resolves a conjecture of Chen and V\'eron. The key step is that the inverse Fourier transform of the high-frequency part of $1/\log(|\xi|^2)$ behaves like $|x|^{-d}(\log(1/|x|))^{-2}$ near the origin.
\end{abstract}

\maketitle

\begingroup
\renewcommand{\thefootnote}{}%
\footnotetext{\emph{Declaration of generative AI and AI-assisted technologies.} During the preparation of this work the authors used Claude (Anthropic) to assist in developing the proof of Theorem~\ref{main}. After using this tool, the authors reviewed and edited the content as needed and take full responsibility for the content of the publication.}%
\endgroup

%=====================================================================
\section{Introduction}\label{sec:intro}
%=====================================================================

The logarithmic Laplacian is a nonlocal operator that arises as the formal derivative of the fractional Laplacian with respect to its order at zero. Namely, as introduced by Huyuan Chen and Tobias Weth~\cite{MR3995092}, it is given formally by
\begin{equation*}
\log(-\Delta)f:=\lim_{\sigma \rightarrow 0^{+}}\frac{1}{\sigma}\left ((-\Delta)^\sigma f-f \right ).
\end{equation*}
One can express the logarithmic Laplacian as
\begin{equation}\label{eqn:logarithmic_laplacian}
\log(-\Delta)f(x)=\gamma_{d}\int_{|x-y|\leq 1}\frac{f(x)-f(y)}{|x-y|^d}\,dy-\gamma_{d}\int_{|x-y|>1}\frac{f(y)}{|x-y|^d}\,dy+\rho_{d} f(x), \quad \text{for $x\in \R^d$,}
\end{equation}
where
\begin{equation*}
\gamma_{d}=\frac{2}{\omega_{d-1}}, \quad \rho_{d}=2\log(2)+\psi(d/2)-\gamma_{E},
\end{equation*}
$\omega_{d-1}$ is the surface measure of the unit sphere $S^{d-1}$, $\psi$ is the digamma function and $\gamma_{E}$ is the Euler--Mascheroni constant. Furthermore, one can justify that the notation $\log(-\Delta)$ is valid due to the identity via the Fourier transform $\F$:
\begin{equation*}
	\F\bigl[\log(-\Delta)f\bigr](\xi)=2\log(|\xi|)\,\F[f](\xi), \quad \text{for $\xi \in \R^d\setminus \{0\}$,}
\end{equation*}
where $f\in \mathcal{S}(\R^d)$ (the space of Schwartz functions). Throughout we use the convention
\begin{equation*}
\F[f](\xi)=\int_{\R^d}f(x)e^{-ix\cdot\xi}\,dx, \qquad \Fi[g](x)=\frac{1}{(2\pi)^d}\int_{\R^d}g(\xi)e^{ix\cdot\xi}\,d\xi.
\end{equation*}

A natural question that arises in the study of the logarithmic Laplacian is the existence of the fundamental solution and its associated properties. The first results in relation to the fundamental solution were established by Chen and V\'eron~\cite{MR4741033}. Namely, they established the existence of the fundamental solution of the logarithmic Laplacian in dimensions $d\geq 3$. In particular, they solved the problem
\begin{equation}\label{eqn:fundamental_soln}
	\log(-\Delta)E=\delta_0 \quad \text{in $\mathcal{D}'(\R^d)$,}
\end{equation}
and showed that the solution is necessarily a locally integrable function.

Here $\mathcal{D}'(\R^d)$ is the space of distributions, which is the dual of the space $\mathcal{D}(\R^d)$ of smooth functions with compact support, and $\delta_0$ is the Dirac measure with mass at the origin.

The result of Chen and V\'eron~\cite[Theorem 1.8]{MR4741033} establishes the existence of the fundamental solution in dimension 3 and higher and gives a rather precise characterization. Naturally, Chen and V\'eron conjectured the following.

\begin{conjecture}[\cite{MR4741033}]\label{conjecture:chen_veron}
	In dimensions $1$ and $2$, the fundamental solution $E$ of the logarithmic Laplacian exists, is a (radial) locally integrable function, and satisfies the bound
	\begin{equation}\label{eqn:conj_bdd}
		|E(x)|\lesssim |x|^{-\frac{d-1}{2}}, \quad \text{for $|x|\geq 2$.}
	\end{equation}
\end{conjecture}

One of the authors of this note was able to establish the existence of the fundamental solution of the logarithmic Laplacian in dimensions 1 and 2~\cite{MR5093271}. However, the notion of solution considered in~\cite{MR5093271} is not exactly the same as the notion of solution in~\eqref{eqn:fundamental_soln}. There \eqref{eqn:fundamental_soln} is solved in the sense of \emph{Lizorkin distributions}, by a division argument in the spirit of the Malgrange--Ehrenpreis theorem, cf.~\cite{MR86990,MR70048,MR68123,MR3995092}. Nevertheless, despite this difference in the notion of solution, one can still effectively identify the fundamental solution as a tempered distribution. In particular, the following was proved.

\begin{theorem}[\cite{MR5093271}]\label{thm:lee}
    There exists a radial tempered distribution $E$ which solves \eqref{eqn:fundamental_soln} in the sense of Lizorkin distributions, such that $E$ can be identified with a function on $\{|x|\geq 2\}$ and $E$ satisfies \eqref{eqn:conj_bdd}.
\end{theorem}

This is almost Conjecture~\ref{conjecture:chen_veron} of Chen and V\'eron. What is missing is that $E$ can be identified with a locally integrable function.

In the proof of~\cite[Theorem 3, Step 2]{MR5093271}, it can be observed that \textbf{to completely resolve the conjecture of Chen and V\'eron (Conjecture~\ref{conjecture:chen_veron})} it is sufficient to verify that the inverse Fourier transform of
\begin{equation*}
    \mathbf{1}_{\{|\xi|>2\}}\frac{1}{\log(|\xi|^2)}, \quad \text{for $\xi \in \R^d$,}
\end{equation*}
is a locally integrable function. This is the content of Theorem~\ref{main}, which allows us to fully verify Conjecture~\ref{conjecture:chen_veron}.

\begin{theorem}\label{main}
Let $\eta \in C^\infty_c(\R^d)$ be a radial cutoff function with $0\le\eta\le 1$ and $\eta \equiv 1$ on $\overline{B_2}$, and let
\begin{equation*}
w := \Fi\left[ (1 - \eta) \frac{1}{\log(|\cdot|^2)} \right].
\end{equation*}
Then $w$ is radial and smooth on $\R^d \setminus \{0\}$, and
\begin{equation*}
w(x) = \frac{\Gamma(d/2)}{4\pi^{d/2}}  \frac{1}{|x|^{d} \left( \log(1/|x|) \right)^{2}} \bigl( 1 + o(1) \bigr), \qquad |x| \to 0.
\end{equation*}
\end{theorem}
The theorem states that the distribution $\Fi\left[ (1 - \eta) \frac{1}{\log(|\cdot|^2)} \right]$ agrees with the function $w$ away from $0$, that is,
\begin{equation*}
\left\langle \Fi\left[ (1 - \eta) \frac{1}{\log(|\cdot|^2)} \right],\varphi \right\rangle = \langle w,\varphi \rangle
\end{equation*}
for all $\varphi\in\mathcal{D}(\R^d)$ with $0\notin\supp(\varphi)$.
\subsection*{Setup}
Let $E$ denote the fundamental solution from Theorem~\ref{thm:lee}, and let $\eta$ be the cutoff function of Theorem~\ref{main}, so that $1 - \eta \equiv 0$ on $B_2$ and $\supp(1 - \eta) \subset \{ |\xi| \geq 2 \}$.

We first remark that
\begin{equation}\label{subord}
\frac{1}{\log(|\xi|^2)} = \int_0^\infty |\xi|^{-2t}\,dt,
\end{equation}
which is valid for all $|\xi| > 1$.

For $0<\alpha<d$ set
\begin{equation*}
k_{d,\alpha} :=  \frac{\Gamma(\frac{d-\alpha}{2})}{\pi^{d/2}\,2^{\alpha}\,\Gamma(\frac{\alpha}{2})}, \qquad k_d(t):=k_{d,2t}=\frac{\Gamma\bigl(\tfrac{d}{2} - t\bigr)}{\pi^{d/2}\, 4^t\, \Gamma(t)}.
\end{equation*}
The \textbf{Riesz potential identity} is
\begin{equation}\label{riesz}
(\mathcal{I}_\alpha \mu)(x) =k_{d,\alpha}\int_{\R^d}\frac{\mu(dy)}{|x - y|^{d-\alpha}},
\end{equation}
where $\mu$ is a compactly supported measure.
By \cite[Ch.~V, \S1.1, Lemma~1]{Stein70},
\begin{equation}\label{steinid}
\Fi\bigl[|\cdot|^{-\alpha}\bigr]=k_{d,\alpha}\,|x|^{\alpha-d}\qquad\text{in }\mathcal{S}'(\R^d),\quad 0<\alpha<d,
\end{equation}
that is,
\begin{equation*}
\bigl\langle k_{d,\alpha}|\cdot|^{\alpha-d},\phi\bigr\rangle=\bigl\langle|\cdot|^{-\alpha},\Fi[\phi]\bigr\rangle \quad \text{for all $\phi\in\mathcal{S}(\R^d)$.}
\end{equation*}

\subsection*{Outline of the proof of Theorem~\ref{main}}
\step{Step 1.} We take the inverse Fourier transform of $(1-\eta)/\log(|\cdot|^2)$ and use \eqref{subord} to split it into two parts,
\begin{equation*}
w_0 = \Fi\left[ (1 - \eta)\int_0^{t_*} |\cdot|^{-2t}\, dt \right], \qquad R = \Fi\left[ (1 - \eta)\frac{|\cdot|^{-2t_*}}{\log(|\cdot|^2)}\right].
\end{equation*}

\step{Step 2.} Use \eqref{steinid} to write
\begin{equation*}
w_0(x) = \int_0^{t_*} \bigl(k_d(t)\,|x|^{2t-d} - r_t(x)\bigr)\,dt, \qquad r_t := \Fi\bigl[\eta\,|\cdot|^{-2t}\bigr],
\end{equation*}
and show that the $r_t$ term is $O(1)$.

\step{Step 3.} We write
\begin{equation*}
\int_0^{t_*} k_d(t)\,|x|^{2t-d}\, dt =|x|^{-d} \int_0^{t_*} k_d(t)\,e^{-tL}\, dt, \qquad L:= 2 \log(1/|x|),
\end{equation*}
and show that
\begin{equation*}
 \int_0^{t_*} k_d(t)\, e^{-tL}\, dt
  = \frac{\Gamma(d/2)}{\pi^{d/2}}\, \frac{1}{L^2} + O\Bigl(\frac{1}{L^3}\Bigr)
  \qquad \text{as } L \to \infty.
\end{equation*}
Thus
\begin{equation*}
w_0(x) = |x|^{-d}\left(\frac{\Gamma(d/2)}{\pi^{d/2}}\, \frac{1}{L^2} + O\Bigl(\frac{1}{L^3}\Bigr)\right) + O(1)
  \qquad \text{as } L \to \infty.
\end{equation*}

\step{Step 4.} We show that the remainder term $R$ is $O(|x|^{2t_*-d})$ and conclude the proof by showing that the $w_0$ term dominates.

\begin{proof}[Proof of Theorem~\ref{main}]\ \\
\vspace{-0.5cm}
\step{Step 1.}
The symbol $(1-\eta)/\log(|\cdot|^2)$ is in $S^0$, so by
\cite[Ch.~VI, \S4.1, Prop.~1]{Stein93}, $w$ is smooth
on $\R^d\setminus\{0\}$ and $|w(x)| \leq A_N |x|^{-d-N}$ for all $N \in \mathbb{N}$, for $x$ away from the origin.

As the support of $1 - \eta$ is contained in $\{ |\xi| \geq 2 \}$, the identity~\eqref{subord} is valid there. Fix $t_* \in (0, d/2)$. Since
\begin{equation*}
\int^{\infty}_{t_{*}}|\xi|^{-2t}\, dt = \frac{|\xi|^{-2t_*}}{\log(|\xi|^2)},
\end{equation*}
splitting the integral in \eqref{subord} at $t_*$ gives
\begin{equation*}
(1 - \eta(\xi)) \frac{1}{\log(|\xi|^2)} = (1 - \eta(\xi))\int_0^{t_*} |\xi|^{-2t}\, dt +  (1 - \eta(\xi))\frac{|\xi|^{-2t_*}}{\log(|\xi|^2)}.
\end{equation*}
We take the inverse Fourier transform of each term and analyse them separately,
\begin{equation*}
w_0 = \Fi\left[ (1 - \eta)\int_0^{t_*} |\cdot|^{-2t}\, dt \right], \qquad R = \Fi\left[ (1 - \eta)\frac{|\cdot|^{-2t_*}}{\log(|\cdot|^2)}\right].
\end{equation*}

\step{Step 2.}
First we justify the interchange of $\Fi$ and the $t$-integral in $\mathcal{S}'$. For $\phi \in \mathcal{S}(\R^d)$ the distributional inverse Fourier transform is
\begin{align*}
\left\langle \Fi\left[(1 -\eta)|\cdot|^{-2t}\right],\phi\right\rangle &= \left\langle (1 -\eta)|\cdot|^{-2t},\Fi[\phi]\right\rangle \\
&= \int_{\R^d}  (1 -\eta(\xi))|\xi|^{-2t}\,\Fi[\phi](\xi)\, d\xi.
\end{align*}
Since $\bigl|(1-\eta(\xi))|\xi|^{-2t}\bigr|  \le 1$ for all $t \in [0,t_*]$, setting $\Phi=\Fi[\phi]\in\mathcal{S}(\R^d)\subset L^1(\R^d)$ we have
\begin{equation*}
\int_0^{t_*} \int_{\R^d} \bigl|(1-\eta(\xi))|\xi|^{-2t}\Phi(\xi)\bigr|\,d\xi\,dt \le t_*\|\Phi\|_{L^1} < \infty,
\end{equation*}
so by Fubini the interchange is justified.

By \eqref{steinid} with $\alpha = 2t$ we have
\begin{equation}\label{3}
\int_0^{t_*}\Fi\bigl[(1-\eta)\,|\cdot|^{-2t}\bigr](x)\, dt=\int_0^{t_*} \bigl(k_d(t)\,|x|^{2t-d} - r_t(x)\bigr)\,dt, \qquad r_t := \Fi\bigl[\eta\,|\cdot|^{-2t}\bigr].
\end{equation}

We now analyse separately the $k_d(t)|x|^{2t-d}$ term and the $r_t(x)$ term.

The function $\eta\,|\cdot|^{-2t}$ has compact support and is in $L^1(\R^d)$ as $2t \leq 2t_* < d$. Thus
\begin{equation*}
|r_t(x)| \le \frac{1}{(2\pi)^d}\int_{\R^d}\bigl|\eta(\xi)|\xi|^{-2t}e^{ix \cdot \xi}\bigr|\,d\xi \leq \frac{1}{(2\pi)^d}\int_{\R^d}\eta(\xi)|\xi|^{-2t}\,d\xi \leq C,
\end{equation*}
with $C$ independent of $t\in[0,t_*]$, and so $\int_0^{t_*}r_t(x)\,dt = O(1)$.

\step{Step 3.}
The first term on the right hand side of \eqref{3} can be written
\begin{equation*}
\int_0^{t_*} k_d(t)\,|x|^{2t-d}\, dt =|x|^{-d} \int_0^{t_*} k_d(t)\,e^{-tL}\, dt, \qquad L:= 2 \log(1/|x|).
\end{equation*}
The function $k_d$ is real analytic near $t= 0$, and since $1/\Gamma(t) = t + O(t^2)$ its Taylor expansion is
\begin{equation*}
k_d(t) = \frac{\Gamma(d/2)}{\pi^{d/2}}\,t + O(t^2), \quad t \to 0^+.
\end{equation*}
Write $k_d(t) = \frac{\Gamma(d/2)}{\pi^{d/2}}\,t + s(t)$. Since
$s(t) = O(t^2)$ as $t \to 0^+$ and $k_d$ is bounded on $[0, t_*]$, there is
a constant $M$ such that $|s(t)| \le M t^2$ for all $t \in [0, t_*]$. Then
\begin{equation*}
  \int_0^{t_*} k_d(t)\, e^{-tL}\, dt
  = \frac{\Gamma(d/2)}{\pi^{d/2}} \int_0^{t_*} t\, e^{-tL}\, dt + \int_0^{t_*} s(t)\, e^{-tL}\, dt .
\end{equation*}
For the first term,
\begin{equation*}
  \int_0^{t_*} t\, e^{-tL}\, dt
  = \frac{1}{L^2} - \Bigl(\frac{t_*}{L} + \frac{1}{L^2}\Bigr) e^{-t_* L}
  = \frac{1}{L^2} + O\bigl(e^{-t_* L}\bigr),
\end{equation*}
and for the second,
\begin{equation*}
  \biggl| \int_0^{t_*} s(t)\, e^{-tL}\, dt \biggr|
  \le M \int_0^{\infty} t^2 e^{-tL}\, dt = \frac{2M}{L^3}.
\end{equation*}
Therefore
\begin{equation*}
  \int_0^{t_*} k_d(t)\, e^{-tL}\, dt
  = \frac{\Gamma(d/2)}{\pi^{d/2}}\, \frac{1}{L^2} + O\Bigl(\frac{1}{L^3}\Bigr)
  \qquad \text{as } L \to \infty.
\end{equation*}

We conclude that
\begin{equation*}
w_0(x) = |x|^{-d}\left(\frac{\Gamma(d/2)}{\pi^{d/2}}\, \frac{1}{L^2} + O\Bigl(\frac{1}{L^3}\Bigr)\right) + O(1)
  \qquad \text{as } L \to \infty.
\end{equation*}

\step{Step 4.}
Finally we analyse the remainder term $R = \Fi\Bigl[ (1 - \eta)\frac{|\cdot|^{-2t_*}}{\log(|\cdot|^2)}\Bigr]$.
The symbol $(1 - \eta(\xi)) \frac{|\xi|^{-2t_*}}{\log(|\xi|^2)}$ is in $S^{-2t_*}$, so by \cite[Ch.~VI, \S4.1, Prop.~1]{Stein93}, $R$ is smooth on $\R^d\setminus\{0\}$ and
\begin{equation*}
|R(x)| \leq A|x|^{2t_*-d},
\end{equation*}
that is, $R(x)=O(|x|^{2t_*-d})$. Since $|x|^{-d}e^{-t_*L}=|x|^{2t_*-d}$, we conclude that
\begin{equation*}
w(x) =  |x|^{-d}\left(\frac{\Gamma(d/2)}{\pi^{d/2}}\, \frac{1}{L^2} + O\Bigl(\frac{1}{L^3}\Bigr)\right) + O(1) + O(|x|^{-d}e^{-t_*L}).
\end{equation*}

It remains to show that the leading term dominates. Let $g(x)=|x|^{-d}L^{-2}$. As $|x|\to 0$ we have $L\to\infty$ and
\begin{equation*}
\frac{|x|^{-d}L^{-3}}{g(x)} = \frac{1}{L} \to 0, \qquad \frac{1}{g(x)} = |x|^{d}L^{2} \to 0, \qquad \frac{|x|^{2t_*-d}}{g(x)} = |x|^{2t_*}L^2 \to 0.
\end{equation*}
Thus, recalling $L = 2\log(1/|x|)$,
\begin{equation*}
w(x)= \frac{\Gamma(d/2)}{\pi^{d/2}}\, |x|^{-d}L^{-2}\bigl(1+o(1)\bigr) = \frac{\Gamma(d/2)}{4\pi^{d/2}}\, \frac{1}{|x|^{d}\left(\log(1/|x|)\right)^{2}}\bigl(1+o(1)\bigr).
\end{equation*}
\end{proof}

\section*{Disclosure on the use of AI}
The proof of Theorem~\ref{main} was developed with the assistance of Claude, an AI model developed by Anthropic. We record two observations on how the argument relates to the existing literature.

First, the focus of~\cite{MR5093271} was the decay of the fundamental solution at infinity, namely the bound~\eqref{eqn:conj_bdd}. That decay is governed by the behaviour of $1/\log(|\xi|^2)$ at low frequencies, and the argument there contains additional steps that serve only this purpose. These steps are not required here. Local integrability concerns the behaviour of $E$ near the origin, which is governed by the high frequencies, and this is why the argument above only takes the inverse Fourier transform of $(1-\eta)/\log(|\cdot|^2)$. Claude did not state this reasoning explicitly. The explanation above is our own, and we suspect it is the reason the argument takes this form.

Second, after comparing the argument with the work of Chen and V\'eron~\cite{MR4741033} on the fundamental solution in dimensions $d\geq 3$, we suspect that Claude largely followed their ideas. Their computations therefore make a substantial contribution to the present note.

\end{document}